\documentclass[11pt,a4paper]{article}
\usepackage{amsmath}
\usepackage{amssymb}
\usepackage{amscd}
\usepackage[latin1]{inputenc}
\usepackage[francais]{babel} 
\usepackage{fancyhdr}

\title{Stabilit\'e de l'in\'egalit\'e de Faber-Krahn \\
en courbure de Ricci positive} \author{J\'er\^ome Bertrand\thanks{Soutenu par la requ\^{e}te 20-101469
du FNRS.}\\
Institut de math\'ematiques, \\ Universit\'e de Neuch\^atel, Suisse. \\
E-mail: jerome.bertrand@unine.ch} \date{\ } 

\newtheorem{lem}{Lemme}[section]
\newtheorem{cor}[lem]{Corollaire}
\newtheorem{theo}[lem]{Th{\'e}or{\`e}me}
\newtheorem{prop}[lem]{Proposition}
\newtheorem{defi}[lem]{D{\'e}finition}
\newtheorem{rem}[lem]{Remarque}

\newenvironment{dem}{\bf Preuve : \rm}{\hfill
$\blacksquare$\par\medbreak} 
\DeclareMathOperator{\hes}{Hess}
\DeclareMathOperator{\vol}{\, \rm{vol} \,}
 \DeclareMathOperator{\im}{\,
\frac{1}{\vol M } \, \int_{M}} 
\DeclareMathOperator{\ens}{\, {\mathcal{M}}_n \,}
\DeclareMathOperator{\ep}{ \varepsilon}
\newcommand{\ov}[1][]{\overline{#1}}

\begin{document}

\maketitle 
\begin{center}
    {\bf Abstract} 
\end{center}
\begin{quotation}
P. Bérard and D. Meyer proved a Faber-Krahn inequality for domains in compact manifolds with positive Ricci curvature. We prove stability results for this inequality.
\end{quotation}
\section*{Introduction} 
En utilisant l'in\'egalit\'e isop\'erim\'etrique de L\'evy-Gromov \cite{gromov99}, P. B\'erard et 
D. Meyer ont démontré une in\'egalit\'e du type Faber-Krahn pour les domaines d'une 
vari\'et\'e compacte à courbure de Ricci positive.
\begin{theo}[Bérard-Meyer, \cite{berard1982}]\label{a4t2}
Soit (M,g) une vari\'et\'e riemannien\-ne compacte de
dimension $n$ dont la courbure de Ricci  v\'erifie
$Ric \geq (n-1)g$. Soit $\Omega$ un domaine r\'egulier de $M$ et
$\Omega^{*}$ le domaine sym\'etris\'e de $\Omega$, c'est-\`a-dire une
boule g\'eod\'esique de la sph\`ere canonique $(\mathbb{S}^n,can)$
v\'erifiant $\frac{\vol (\Omega)}{\vol (M)}=\frac{\vol
(\Omega^{*})}{\vol (\mathbb{S}^n)}$ . Sous ces hypothèses, on a l'in\'egalit\'e 
 $$
\lambda_1^D (\Omega) \geq \lambda_1^D(\Omega^{*}), $$ 
o\`u $ \lambda_1^D$ désigne la première valeur propre de Dirichlet du domaine
sur l'espace correspondant. De plus, l'\'egalit\'e a lieu si et
seulement si le triplet $(\Omega,M,g)$ est isom\'etrique au triplet ($
\Omega^{*},\mathbb{S}^n,can$). 
\end{theo} 

L'objet de cet article est d'\'etudier les domaines des vari\'et\'es \`a courbure 
de Ricci positive dont la premi\`ere valeur propre de Dirichlet est proche de celle 
de leur domaine symétrisé.

La premi\`ere remarque est que de tels domaines ne sont pas n\'ecessairement 
hom\'eomorphes \`a des boules euclidiennes. En effet si l'on retire des ensembles de petite
capacit\'e \`a une calotte sph\'erique de la sph\`ere canonique (par exemple des boules de petit rayon), on modifie peu la premi\`ere valeur propre de  Dirichlet \cite{chavel78,rauch75}. Il est  \'egalement facile de construire des exemples de vari\'et\'e qui ne sont pas proches,
pour la distance de Gromov-Hausdorff, de la sph\`ere canonique et qui pourtant contiennent
des domaines dont la premi\`ere valeur propre de Dirichlet est arbitrairement 
proche de celle de leur domaine sym\'etris\'e (par exemple en lissant aux extr\'emit\'es 
des sinus produits tordus de la forme $((0,\pi)\times \mathbb{S}^{n-1}, dt^{2}+\ep^2\sin^{2}t \,can)$ avec $\ep>0$ petit et en considérant des domaines de la forme $(0,a) \times \mathbb{S}^{n-1}$). Signalons enfin qu'il est \'egalement
possible de construire de tels exemples sur des vari\'et\'es \`a courbure de Ricci positive,
non hom\'eomorphes \`a la sph\`ere en utilisant les m\'etriques construites sur l'espace
projectif complexe par M. Anderson \cite{anderson1990}. Muni de ces métriques l'espace projectif complexe est proche, pour la distance de Gromov-Hausdorff, d'un sinus produit tordu, pour plus de détails sur cet exemple nous renvoyons à \cite{these}.
   
 Dans cet article nous démontrons un résultat de stabilité optimal, au vue des remarques précédentes, 
 lorsque la première valeur propre d'un domaine convexe est proche de celle d'un hémisphère $\mathbb{S}^{n+}$ de la sphère canonique de dimension $n$.  
 
\begin{theo}\label{a8t1}Il existe des fonctions $\eta(\ep)$ et $\tau(\ep)$ 
telles que, pour toute vari\'et\'e riemannienne compacte $(M^n,g)$ dont la courbure de Ricci vérifie $Ric  \geq (n-1)g$, pour tout domaine régulier $\Omega$ de $M$, géodésiquement convexe, de
volume $ \vol \Omega \leq \frac{1}{2} \vol M $ et dont la premi\`ere
valeur propre de Dirichlet v\'erifie
$$\lambda_1^D(\Omega) \leq \lambda_1^D(\mathbb{S}^{n+}) + \epsilon, $$ 
alors, en notant $d_{H}$ la distance de Hausdorff, il existe $x_{0}$ dans $\Omega $ tel 
que 
$$ d_{H} (\Omega, B(x_0, \frac{\pi}{2})) \leq \tau (\ep).$$
De plus, il existe un sinus produit tordu $((0,\pi)\times N, d)$ tel que
$$\begin{array}{l}
 d_{GH}((M,d_{g}),((0,\pi)\times N, d)) \leq \tau(\ep),\\
 d_{GH} ((\Omega, d_g),
 ((0,\frac{\pi}{2})\times N, d)) \leq \tau(\ep),
\end{array}$$
où $d_{g}$ est la distance induite par la m\'etrique riemannienne $g$ et $d_{GH}$ désigne la distance de Gromov-Hausdorff. Nous renvoyons à la  troisième partie de cet article pour une définition des sinus produits tordus.
\end{theo} 

Dans la deuxième partie de cet article, nous démontrons également un résultat de stabilité plus faible mais sous une hypothèse de courbure moyenne du bord positive (théorème \ref{a3t2}).

 Dans le théorème \ref{a8t1}, l'hypothèse sur le volume relatif du domaine implique par le résultat de P. Bérard et D. Meyer que la première valeur propre est supérieure ou égale à celle d'un hémisphère. Cet hypothèse sur le volume relatif est nécessaire pour obtenir le résultat sur la variété ambiante. Pour s'en convaincre, il suffit de considérer un produit tordu de la forme $((0,\frac{\pi}{2})\times \mathbb{S}^{n-1}, dt^2+\ep^2\sin^2 t can)$ que l'on recolle avec un hémisphère de $(\mathbb{S}^{n}, \ep^2 can)$. 
 \smallskip
Le th\'eor\`eme \ref{a8t1} g\'en\'eralise dans le cas de l'h\'emisph\`ere, un r\'esultat de A. Avila
\cite{avila2002} sur certains domaines convexes de la sph\`ere $\mathbb{S}^2$.

\begin{theo}[Avila]Soit $\Omega$ un domaine r\'egulier géodésiquement con\-vexe contenu dans un h\'emisph\`ere de $\mathbb{S}^2$. Soit $B$ une boule de m\^eme volume que $\Omega$. Supposons que 
$$ \lambda_1^D(\Omega) \leq \lambda_1^D(B) + \ep \, ,$$
alors il existe une fonction $\tau(\ep)$ d\'ependant de $\vol(\Omega),\vol (\partial \Omega)$ et du rayon  de $B$ telle que  $\Omega$ est $\tau(\ep)$-Hausdorff proche de $B$.
\end{theo}

Un des éléments de la preuve du théorème \ref{a8t1} est de montrer que la première fonction propre du domaine considéré est proche de la fonction propre correspondante dans le cas modèle. Nous avons regroupé ces résultats dans la première partie de cet article. Dans la deuxième partie, nous démon\-trons un résultat de stabilité plus faible que le théorème \ref{a8t1} mais où l'hypothèse sur la convexi\-té est remplacée par une hypothèse (plus faible) de courbure moyenne du bord positive (théorème \ref{a3t2}). Nous utilisons ensuite ce résultat pour démon\-trer dans une troisième partie le théorème \ref{a8t1}.

\section{Résultats pr\'eliminaires}
\begin{defi}On note $\ens$ l'ensemble (des classes d'isométrie) des vari\'e\-t\'es riemanniennes connexes, compactes, 
de dimension $n$ dont la courbure de Ricci v\'erifie $Ric \geq (n-1)g$.
\end{defi}

Soit $ p \geq 1$ un nombre r\'eel et $h$ appartenant à $L^p (M)$. On note
$$ ||h||_{L^p}= \left(\frac{1}{\vol M} \int_M |h|^p \,
dx\right)^{\frac{1}{p}}. $$
On utilisera la définition usuelle pour la norme $L^{\infty}$.

On notera $\tau (\ep),r(\ep), \eta(\ep)$, etc $\dots$ de mani{\`e}re 
g{\'e}n{\'e}rique, toute \\quantit{\'e} positive ne
d{\'e}pendant que de $\ep$ et de la dimension $n$ de la 
vari{\'e}t{\'e}, dont la limite quand $\ep$ tend vers $0$ est $0$.

Dans la suite, on suppose que la premi\`ere fonction propre de Dirichlet sur un domaine
$\Omega$ est normalis\'ee par 
\begin{equation}\label{a4e4}
    \sup_{\Omega}f= 1.
\end{equation}

\subsection{Formule de Reilly}
Les fonctions propres de la sph\`ere canonique $\mathbb{S}^n$ de valeurs propres
 $\lambda_1(\mathbb{S}^n) \\=n$ v\'erifient l'\'equation
\begin{equation}\label{obata}
 \hes f + fcan =0 
\end{equation}
 avec $can$ la m\'etrique canonique de la sph\`ere. Il en est de m\^eme pour la premi\`ere fonction propre de Dirichlet d'un 
 h\'emisph\`ere. C'est une cons\'equence du lemme suivant (nous renvoyons à \cite{chavel84} pour 
 plus de détails).
 
 \begin{lem}Soit $f$ une fonction propre associ\'ee \`a la premi\`ere valeur propre non nulle
     $\lambda_{1}(M)$ d'une vari\'et\'e riemannienne compacte $(M,g)$. Soit $\Omega$
     un domaine nodal de la fonction $f$ alors
     $$ \lambda_{1}^{D}(\Omega)= \lambda_1 (M)$$
et la premi\`ere fonction propre du domaine est la restriction de $f$ \`a $\Omega$.
 \end{lem}
\begin{rem}On peut définir la première valeur propre de Dirichlet même lorsque le bord 
du domaine n'est pas régulier.
\end{rem}  
Sous les hypothèses du théorème \ref{a8t1}, la norme $L^2$ du membre de gauche de l'équation (\ref{obata}) reste petite. Cette estimation sur le hessien de la première fonction propre est une conséquence d'une formule due \`a R. Reilly \cite{reilly}.

\begin{lem}[Formule de Reilly]
Soit $(M^n,g)$ une vari\'et\'e riemannienne compacte éventuellement \`a bord
lisse $\partial M$. Pour toute fonction $f$ appartenant \`a
$C^{\infty}(\ov[M])$, on a 
\begin{multline} \int_M |\hes f |^2-
\int_M (\Delta f)^2 + \int_M \rm{Ric}(\nabla f, \nabla f) \\ = -2
\int_{\partial M}<\nabla^{\partial M}(\frac{\partial f}{\partial
\nu}),\nabla^{\partial M}f> - \int_{\partial M} H \left(\frac{\partial
f}{\partial \nu}\right)^2- \int_{\partial M} \Pi(\nabla^{\partial
M}f,\nabla^{\partial M}f) 
\end{multline} 
o\`u $\ov[M]= M \cup \partial M$, $\nabla^{\partial M}$ d\'esigne le gradient pour la m\'etrique 
induite par $g$ sur $ \partial M$, où la courbure moyenne $H$ est la trace de la seconde forme fondamentale définie par
$\Pi(X,Y)= -g(D_X \eta,Y)$ avec $\eta$ la normale unitaire rentrante et où le membre de droite
 est nul si la variété est sans bord.
 \end{lem}
 
     On d\'eduit de la formule de Reilly le
    
\begin{lem}\label{a4l1}Soit $(M^{n},g)$ une vari\'et\'e riemannienne compacte 
 (respectivement \`a bord, dont la courbure moyenne du bord est positive ou nulle) dont la courbure de Ricci vérifie $Ric \geq (n-1)g$. Notons $\lambda$ la premi\`ere valeur propre non nulle 
 (respectivement  de Dirichlet) et soit $f$ une fonction propre associ\'ee \`a 
 cette valeur propre. Supposons que cette valeur propre v\'erifie
 $$ n \leq \lambda \leq n+ \ep $$
 pour $0<\ep<1$, alors il existe une constante $C(n)$ telle que
 $$ || \hes f + fg||_{L^{2}} \leq C(n) \ep^{\frac{1}{2}}||f||_{L^{2}}.$$
 \end{lem}

 \begin{dem} 
La formule de Reilly appliqu\'ee \`a $f$
donne  
 $$\int_M |\hes f |^2- \int_M (\Delta f)^2 +
\int_M \rm{Ric}(\nabla f, \nabla f) = - \int_{\partial M} H
\left(\frac{\partial f}{\partial \nu}\right)^2.$$ En utilisant
l'hypoth\`ese sur la courbure de Ricci et sur la courbure moyenne du
bord $\partial M$, il vient  
$$ \int_M |\hes f |^2
+\left((n-1)\lambda- \lambda^2\right) \int_M f^2   \leq 0.
$$
On écrit ensuite le terme $\hes f$ sous la forme  $
\hes f= (\hes f +\frac{\lambda }{n} fg) - \frac{\lambda }{n} fg$. Le premier terme {\'e}tant de trace nulle, les deux termes sont orthogonaux pour le produit scalaire usuel sur $L^2(M)$, l'hypothèse sur $\lambda$ permet alors de conclure.
 
\end{dem}

Un r\'esultat du \`a J. Cheeger et T. Colding (\cite{cheeger95}, théorème 2.11) 
permet de d\'eduire des informations g\'eom\'etriques de cette in\'egalit\'e sur
le hessien.

\begin{lem}[\cite{cheeger95}]\label{a1l1}Soit $ (M,g)$ un élément de $\ens$. Il existe des constantes ne d{\'e}pen\-dant que de $n$ not\'ees $
  C(n) $ et $ \tilde{C}(n) $ telles que pour tout ouvert $U_{1}$ et $U_{2}$ de $M$ et 
  pour toute fonction continue $f$ sur $M$, on a 
  
\begin{multline}\label{tt}
 \frac{1}{\vol(U_1 \times U_2)} \int_{ U_1 \times U_2}
\left(  \int_{0}^{1} |(f\circ \gamma_{xy})''(t) +(f\circ \gamma_{xy})(t)|^2dt\right)
dxdy  \\ 
\leq C(n) \left( \frac{1}{\vol U_1} + \frac{1}{\vol U_2} \right)
\int_{M}|\hes(f)+fg|^{2}.
\end{multline}
On obtient dans le cas particulier o\`u $U_{1}$ et $U_{2}$ sont deux boules 
g\'eod\'esiques de rayon $r$ 
\begin{multline}
 \frac{1}{\vol(U_1 \times U_2)} \int_{ U_1 \times U_2}
 \left(\int_{0}^{1} |(f\circ \gamma_{xy})''(t) +f\circ \gamma_{xy}(t)|^2dt
 \right)dxdy \\
 \leq \frac{\tilde{C}(n)}{V(r)} ||\hes(f)+fg||_{L^{2}}^{2}, 
 \end{multline}
où $V(r)$ d\'esigne le volume d'une boule 
g\'eod\'esique de $(\mathbb{S}^n,can)$ de rayon $r$.
\end{lem}

\begin{rem} La notation $ U_1 
\times U_2$ 
d{\'e}signe en r{\'e}alit{\'e} le
sous-ensemble de mesure pleine de ce produit, constitu{\'e} par les
couples $(x,y)$ admettant une unique g{\'e}od{\'e}sique minimisante les 
reliant (not\'ee $ {\gamma}_{xy}$). Soit $W$ un ouvert de $M$ tel que toute géodésique minimisante dont les extrémités appartiennent à $U_1\times U_2$ est contenue dans $W$. On peut remplacer dans (\ref{tt}),
$\int_{M}|\hes(f)+fg|^{2}$ par $\int_{W}|\hes(f)+fg|^{2}$, c'est ce que nous ferons lorsque nous utiliserons ce résultat pour une fonction propre sur un domaine $\Omega$ d'un élément de $\ens$.
\end{rem}
\begin{dem}
Ce r\'esultat est une application directe du th\'eor\`eme 2.11 de \cite{cheeger95} \`a 
la fonction $|\hes f + f g|^{2}$, en remarquant  que pour toute g\'eod\'esique 
$\gamma$ param\'etr\'ee par longueur d'arc, on a
$$ |(f\circ \gamma)''(t) +(f\circ \gamma) (t)|^2 \leq  |\hes f 
+fg|^2(\gamma (t)). $$
\end{dem}
En appliquant le lemme \ref{a1l1} à une fonction $f$ vérifiant les hypothèses du lemme \ref{a4l1} pour des boules $B_1$ et $B_2$ de rayon $r(\ep)$ convenable (supposons que $r(\ep)$ vérifie $\frac{\ep}{V(r(\ep))} \leq \ep^\frac{1}{2}$ et que $||f||_{L^2} \leq 1$), on déduit du lemme \ref{a4l1}  l'in\'egalit\'e
\begin{equation}\label{rajout1}
  \mbox{$\frac{1}{\vol(B_1 \times B_2)}$} \int_{ B_1 \times B_2}
 \int_0^{d(x,y)} |(f\circ \gamma_{xy})''(t) +f\circ 
\gamma_{xy}(t)|^2dtdxdy
 \leq 
\tilde{C}(n) \ep^\frac{1}{2}.
\end{equation}
Notons 
$$C=\left\{ (x,y) \in B_1 \times B_2; \;
 \int_0^{d(x,y)} |(f\circ \gamma_{xy})''(t) +f\circ 
\gamma_{xy}(t)|^2dt \leq \ep^\frac{1}{4}\right\}.$$
 On déduit de (\ref{rajout1}) l'estimation (inégalité de Byenaim\'e-Tchebitchev) 
 $$ \vol (C) \geq (1- \tilde{C}(n)\ep^\frac{1}{4})\vol (B_1\times B_2).$$
Pour $\ep$ assez petit, l'ensemble $C$ est donc non vide et on en déduit l'existence de couples $(x,y) $ appartenant à $B_1 \times B_2$ pour 
lesquels $ f\circ \gamma_{xy}$ v\'erifie presque la m\^eme \'equation 
diff\'erentielle que dans le cas de la sph\`ere canonique. On peut ensuite par des m\'ethodes classiques comparer $ f\circ 
\gamma_{xy}$ \`a une solution correspondante sur la sph\`ere en 
fixant des conditions au bord \`a l'aide du 
lemme suivant.

 \begin{lem}\label{a1l9.5}Soit $v(t)$ et $Z(t)$ deux fonctions 
d{\'e}finies
  sur $ [0,l]$ avec $l<\pi$.  On suppose que $ \int_0^l Z^2(t)dt < 
\epsilon^2 $ et que
 $v$ est solution de $ v"+v=Z $ avec $ |v(0)-a|< \eta $ et   $
 |v(l)-b| < \eta $. Il existe une constante positive $C$ telle que pour tout $t$ dans $[0,l]  $,
 $$ |v(t)-\tilde{u}_{a,b}(t)| < \frac{C}{\sin (l)}(\ep +\eta)  $$
 et
 $$|v'(t)-\tilde{u}_{a,b}'(t)| <\frac{C}{\sin (l)}(\ep +\eta),$$
o\`u $\tilde{u}_{a,b}  $ est la solution de $ {u}" + u=0$ sur $[0,l] 
$ v{\'e}rifiant les conditions initiales $ u(0)=a $ et $ u(l)=b $.
\end{lem}
On peut \'egalement fixer des conditions de Cauchy (pour une démonstration des lemmes \ref{a1l9.5} et \ref{a1l9}, nous renvoyons à \cite{these}).
\begin{lem}\label{a1l9}Soit $v(t)$ et $Z(t)$ deux fonctions 
d{\'e}finies
  sur $ [0,l]$ avec $l \leq \pi$. On suppose que $ \int_0^l Z^2(t)dt 
< \epsilon^2 $ et que
 $v$ est solution de $ v"+v=Z $ avec $ |v(0)-a|< \eta $ et  $ 
|v'(0)-b|<\eta $
 . Il existe une constante positive $C$ telle que pour tout $t$ dans $[0,l], $
 $$ |v(t)-u_{a,b}(t)| < C(\epsilon +\eta)$$
 et 
$$|v'(t)-u_{a,b}'(t)| < C(\epsilon+\eta) ,$$
o\`u $u_{a,b}$ est la solution de $ {u}" + u=0$ sur $[0,l] $ 
v{\'e}rifiant les conditions initiales $ u(0)=a $ et
$ u'(0)=b$.
\end{lem} 

Pour contr\^oler les conditions initiales de l'\'equation 
diff\'erentielle dans le lemme \ref{a1l9}, nous aurons 
besoin d'une estimation qui prouve que la norme du gradient d'une fonction 
propre sur les domaines consid\'er\'es, reste petite au voisinage
des points r\'ealisant les extr\'ema de la fonction propre. Cette estimation 
est du type de celles obtenues par P. Li et S.T. Yau \cite{li1980}.

\begin{prop}\label{a3l5}Soit $(M,g)$ un élément de $\ens$ et $\Omega$ un domaine r\'egulier de M dont la
courbure moyenne en tout point du bord $\partial \Omega$ est positive
ou nulle. Soit f la premi\`ere fonction propre de Dirichlet sur
$\Omega$, que l'on suppose normalis\'ee par (\ref{a4e4}). 
Alors pour tout $x$ dans $\Omega$,
  $$ |\nabla f|^2(x) \leq
\lambda_1^D(\Omega)(1-f^2(x)),$$ 
en particulier 
 $$ || \nabla f
||_{L^{\infty}}\leq \sqrt{\lambda_1^D(\Omega)}.$$
 \end{prop}

\begin{dem} Commen\c{c}ons par quelques remarques sur $f$. On note
$\eta (x) $ la normale intérieure unitaire en $ x$ apppartenant à $\partial
\Omega$. En appliquant le principe du maximum fort \`a $-f$, on en
d\'eduit pour tout $x$ dans $\partial \Omega$, 
 $$ \frac{\partial f}{\partial \eta} (x) >0.$$ 
Par cons\'equent en tout point $x$ dans $\partial \Omega$ 
$$ \eta(x)=\frac{\nabla f}{|\nabla f|}(x).$$ 
Donc pour $X,Y$ appartenant à $T_x \, \partial \Omega$ 
 \\ $$\Pi(X,Y)=
\frac{-1}{|\nabla f|(x)}\hes f (X,Y).$$
On en d\'eduit
l'expression suivante de la courbure moyenne, 
 $$H(x)=\frac{-1}{|\nabla
f|(x)}\sum_{i=1}^{n-1}\hes f(e_i,e_i),$$ 
avec $(e_i)_{1 \leq i\leq
n-1}$ une base orthonorm\'ee de $T_x\partial \Omega$. Par
continuit\'e, pour $x$ appartenant à $\partial \Omega$,  $\Delta f(x)=\lambda_1^D(\Omega)f(x)=0$, d'o\`u 
$$H(x)= \frac{1}{|\nabla f|(x)}\hes
f(\eta(x),\eta(x)).$$ 
En particulier, pour tout $x$ dans $\partial \Omega$
\begin{equation}\label{a3e11} \hes f(\eta(x),\eta(x))\geq 0.
\end{equation} 
Passons maintenant \`a la d\'emonstration du lemme. On introduit la fonction 
 $$ F=\frac{|\nabla f|^2}{\beta-f^2}$$
avec $\beta > 1$ un réel fixé. Par compacit\'e, il existe $x_0$ tel que
$F(x_0)=\sup_{\ov[\Omega]}F$. Supposons tout d'abord que $x_0 $
appartient \`a $\partial \Omega$. Dans ce cas,
par le principe du maximum, on doit avoir l'estimation 
 $$\frac{\partial F}{\partial \eta}(x_0) \leq 0.$$

Calculons $ \frac{\partial F}{\partial \eta}(x_0)$.
$$ \frac{\partial F}{\partial \eta}(y)=
\frac{2<D_{\eta} \nabla f,\nabla f>
}{\beta-f^2}+ \frac{2f|\nabla f|^2 \frac{\partial f}{\partial
\eta}(y)}{(\beta-f^2)^2}. $$ 
Or, au point $x_0$
 $$\frac{2f|\nabla f|^2 \frac{\partial
f}{\partial \eta}}{(\beta-f^2)^2} = 0,$$
 puisque $f$ v\'erifie les conditions de Dirichlet sur le bord.\\
Montrons que le premier terme est positif ou nul en $x_0$ (dans ce qui suit, $\eta$ désigne $\eta(x_0)$).
 $$
<D_{\eta} \nabla f,\nabla f>= \hes f(
\eta, \nabla f)$$
$$<D_{\eta} \nabla f,\nabla f>= |\nabla
f|\hes f(  \eta,\eta).$$ 
Donc par (\ref{a3e11}), 
 $$<D_{\eta} \nabla f,\nabla f>(x_0) \geq 0$$ 
et par cons\'equent
\begin{equation}\label{a8e5} \frac{\partial F}{\partial \eta}(x_0)=0.
\end{equation}
 On d\'eduit de (\ref{a8e5}) et du principe du maximum fort l'in\'egalit\'e 
  $$ \Delta
F(x_0) \geq 0.$$

Calculons maintenant le laplacien de $F$.
 $$ dF = \frac{d(|\nabla
f|^2)}{\beta -f^2} + \frac{2f |\nabla f|^2 df}{(\beta-f^2)^2}.$$ 
D'o\`u
 \begin{multline*} \hes F = \frac{\hes (|\nabla f|^2)}{\beta -f^2} +
\frac{4f }{(\beta-f^2)^2}d(|\nabla f|^2)\otimes df \\ + \frac{2 |\nabla
f|^2}{(\beta-f^2)^2}df \otimes df + \frac{8f^2 |\nabla f|^2}{(\beta
-f^2)^3}df \otimes df + \frac{2f |\nabla f|^2}{(\beta-f^2)^2}\hes f.
\end{multline*} 
On en d\'eduit 
 \begin{multline*} \Delta F = \frac{ \Delta( |\nabla
f|^2)}{\beta-f^2} -\frac{4f}{(\beta-f^2)^2}g(d(|\nabla f|^2),df)
-\frac{2 |\nabla f|^4}{(\beta-f^2)^2}- \frac{8f^2 |\nabla f|^4}{(\beta
-f^2)^3} \\ + \frac{2f |\nabla f|^2 \Delta f}{(\beta -f^2)^2}.
\end{multline*}
En particulier, en un point $y$ appartenant à $\partial \Omega$ 
\begin{equation}\label{a4e9} 
 \Delta F(y) = 
\frac{ \Delta( |\nabla f|^2)(y)}{\beta}-\frac{2 |\nabla
f(y)|^4}{\beta^2}.
\end{equation}
 Or, d'apr\`es la formule de Bochner  $$
\frac{1}{2}\Delta (|\nabla f|^2) = -|\hes f|^2 -Ric(\nabla f,\nabla f)
+ \lambda_1^D(\Omega)|\nabla f|^2,$$ puisque $f$ est une fonction
propre de valeur propre $\lambda_1^D(\Omega)$. Par hypothèse sur la courbure, on en déduit 
 $$\Delta (|\nabla f|^2) \leq 2 \lambda_1^D(\Omega) |\nabla f|^2,$$ 
 que l'on injecte dans (\ref{a4e9}) pour obtenir
  \begin{equation}\label{a8e6} \Delta F (y) \leq 2 \left(
\lambda_1^D(\Omega) \frac{|\nabla f|^2}{\beta} - \frac{|\nabla
f|^4}{\beta^2}\right). \end{equation} 
En appliquant (\ref{a8e6}) au
point $x_0$, on en d\'eduit puisque $F(x_0)\geq 0$, 
 $$ F(x_0) \leq
\lambda_1^D(\Omega),$$ 
d'o\`u l'estimation en faisant tendre $\beta$
vers 1.

Nous renvoyons \`a \cite{li1980} pour le cas où $x_{0}$ appartient à 
$\Omega$.
\end{dem}

\section{Domaines \`a courbure moyenne positive}

\subsection{Approximation de Hausdorff et sinus produit tordu}
Pour estimer la distance de Gromov-Hausdorff entre deux espaces m\'etri\-ques, nous 
utiliserons des $\ep$-approximations de Hausdorff dont nous rappelons la d\'efinition ci-dessous.

\begin{defi}[$\ep$-approximation de Hausdorff]
Soit $(X,d)$ et $(Y,\delta)$ \\deux espaces m\'etriques compacts. Une 
$\ep$-approximation de Hausdorff de $X$ dans $Y$ est une application non
 n\'ecessairement continue $ \phi : X \rightarrow Y$ telle que
 $$\phi(X) \mbox{ est un } \ep \mbox{ r\'eseau de } Y$$
 et pour tout $x,x'$ dans $X$
  $$ | d(x,x')- \delta(\phi(x),\phi(x'))| \leq \ep.$$
\end{defi}
Lorsqu'il existe une $\ep$-approximation de Hausdorff entre deux espaces métri\-ques, la distance de Gromov-Hausdorff entre ces deux espaces est majorée par $5\ep$ \cite{gromov99}.

Les espaces mod\`eles apparaissant dans nos r\'esultats de stabilit\'e sont des espaces
de longueur model\'es à partir de la sph\`ere canonique (privée de deux points antipodaux) décrite en coordonnées géodésiques par rapport à un point, appel\'es sinus produit tordu dont nous 
rappelons également la définition.
\begin{defi} Soit $(N,\delta)$ un espace de longueur de diamètre inférieur à $\pi$, on appelle sinus produit tordu, 
l'espace de longueur $((0,\pi)\times N, d)$ o\`u la distance $d$ est d\'efinie 
pour $(t,x), (s,y)$ dans $(0,\pi)\times N$, par
\begin{equation}\label{a4e10}
\cos d((t,x),(s,y))= \cos s \cos t + \sin s \sin t \cos \delta (x,y).
\end{equation}
\end{defi}

Par d\'efinition de la distance, les sous-ensembles de la forme $(0,a)\times N$ avec $a \leq \frac{\pi}
{2}$ forment des parties convexes de  $((0,\pi)\times N, d)$. Par analogie avec le cas de 
la sph\`ere, on appellera h\'emisphère d'un sinus produit tordu la partie $(0,\frac{\pi}
{2})\times N$. Cette propri\'et\'e de  convexit\'e est presque conserv\'ee dans les 
r\'esultats de stabilit\'e que nous avons obtenus. Pour pr\'eciser cette notion 
de presque convexit\'e, nous avons besoin de la d\'efinition qui suit.  

\begin{defi}Soit $(M,g)$ une vari\'et\'e riemannienne compacte. Soit $x$ appartenant à $M$ et $r$ 
un r\'eel positif. On note $d^{x,r}$ la distance induite par la m\'etrique riemannienne sur
la boule g\'eod\'esique $B(x,r)$.
\end{defi}

\subsection{\'Enonc\'e des r\'esultats}
Dans cette partie, nous démontrons le 

\begin{theo}\label{a3t2}Il existe des fonctions $\eta(\ep)$ et $\tau(\ep)$ 
telles que, pour tout élément $(M,g)$ de $\ens$ et pour tout domaine régulier $\Omega$ de $M$, dont
la courbure moyenne $H$ est positive ou nulle en tout point du bord, dont le 
volume vérifie $ \vol \Omega \leq \frac{1}{2} \vol M $ et dont la premi\`ere
valeur propre de Dirichlet v\'erifie
$$\lambda_1^D(\Omega) \leq n+\epsilon, $$ 
alors il existe $x_{0}$ dans $\Omega $ tel que $\Omega$ contient $B(x_{0},\frac{\pi}{2}-\eta(\ep))$ et
$$ \vol ( \Omega \setminus B(x_0,\pi/2-\eta (\ep))) \leq \tau (\ep).$$
De plus, il existe un sinus produit tordu $((0,\pi)\times N, d)$ tel que
$$ d_{GH}((M,d_{g}),((0,\pi)\times N, d)) \leq \tau(\ep)$$
(avec $d_{g}$ la distance induite par la m\'etrique riemannienne) et
$$ d_{GH}((B(x_{0},\frac{\pi}{2}-\eta(\ep)),d^{x_{0},\frac{\pi-\eta(\ep)}{2}}),
 ((0,\frac{\pi}{2})\times N, d)) \leq \tau(\ep).
$$
\end{theo}
\subsection{D\'emonstration du th\'eor\`eme \ref{a3t2}}
Dans le prochain paragraphe, nous démontrons que le domaine $\Omega$ contient une boule de rayon proche de $\frac{\pi}{2}$. Dans le paragraphe suivant, nous montrons à l'aide d'un résultat de J. Cheeger et T. Colding, qu'il existe une approximation de Hausdorff de la variété ambiante dans un sinus produit tordu et que l'image de $\Omega$ par cette approximation, contient une partie qui est Hausdorff proche d'un hémisphère du produit tordu. Dans le dernier paragraphe, nous démontrons une propriété de presque convexité (\ref{a4e5}). Dans la suite, on note $f$ la premi\`ere fonction propre de Dirichlet sur le  domaine $\Omega$,
normalis\'ee par (\ref{a4e4}) et $x_0$ un point de $\Omega$  tel que $f(x_0)=1$.

\subsubsection{$\Omega$ contient une boule de rayon presque \'egal \`a $\frac{\pi}{2}$}
Soit $\Omega$ un domaine r\'egulier de
$(M,g)$ v\'erifiant les hypoth\`eses du th\'eor\`eme \ref{a3t2}.
Notons $d$ la distance de $x_0$ au bord $\partial \Omega$. Par définition, la boule $ B(x_0,d)$ est contenue dans $\Omega.$ L'hypoth\`ese $\frac{\vol \Omega}{\vol M}\leq \frac{1}{2}$ implique par le th\'eor\`eme de Bishop-Gromov, $d\leq \frac{\pi}{2}.$ 
\begin{lem}\label{a3l2}Soit $\Omega$ un domaine
r\'egulier de $(M,g)$ v\'erifiant les hypoth\`eses du th\'eor\`eme
\ref{a3t2}. 
Alors il existe
une fonction $\eta(\ep)$, telle que
pour tout domaine $\Omega$ et tout $(M,g)$ vérifiant les hypothèses
ci-dessus, on a 
 $$
B\left(x_0,\frac{\pi}{2}- \eta(\ep)\right) \subset \Omega.$$
\end{lem}

\begin{dem} 
Sous ces hypoth\`eses, on a par le lemme \ref{a4l1}
  $$||\hes f + fg||_{L^2} \leq C(n) \ep^{\frac{1}{2}}.$$
La première étape consiste \`a appliquer le lemme \ref{a1l1} sur des boules telles,
 qu'une g\'eod\'esique minimisante ayant ces extr\'emit\'es dans celles-ci, soit 
 contenue dans $\Omega$. Soit $y$ appartenant à $\partial \Omega$ tel que $d(x_0,y)=d(x_0,\partial\,\Omega)$,
 $\gamma$ une g\'eod\'esique minimisante reliant $x_0$ \`a $y$ et  $r(\ep)$ un
r\'eel positif tel que  $\frac{\ep}{V(r(\ep))} \leq \ep^{\frac{1}{2}}$ 
($V(r)$ est le volume d'une calotte sph\'erique de rayon $r$ de $\mathbb{S}^{n}$).
Consid\'erons $z$ appartenant à $\gamma$ tel que $d(y,z)=3 r(\ep)$ et $B_1=B(x_0,r(\ep))$, $B_2=B(z,r(\ep)).$ Toute g\'eod\'esique minimisante ayant ces extr\'emit\'es dans $B_1$ et
$B_2$ est n\'eces\-sairement contenue dans $\Omega$. En effet, par l'in\'egalit\'e triangulaire
et pour tout $(u,v)$ dans $B_{1}\times B_{2}$, on a
 $$ d(u,\partial \Omega) >
d-r(\ep) \mbox{ et }
  d(u,v)\leq  d-r(\ep).$$
Le lemme \ref{a1l1} appliqu\'e \`a $B_{1}$ et $B_{2}$ donne 
$$ \frac{1}{\vol(B_1 \times B_2)} \int_{ B_1 \times B_2}
 \left(\int_{0}^{1} |(f\circ \gamma_{uv})''(t) +f (\gamma_{uv}(t))|^2dt
 \right)dudv
 \leq \tilde{C}(n)\ep^{\frac{1}{2}}.$$
 L'\'etape suivante consiste \`a appliquer le lemme \ref{a1l9}. Pour cela on estime
 les conditions initiales v\'erifi\'ees par $f\circ \gamma_{uv}$. Par choix de la normalisation, $f(x_0)=1$ et par la proposition \ref{a3l5}, pour tout $x$ dans $\Omega$,
 \begin{equation}\label{equath2.4}
    |\nabla f|^2(x) \leq
(n+\ep)(1-f^2(x)),
\end{equation}
donc, pour tout $u$ dans $B_1$ et pour tout $v$, $( f\circ \gamma_{uv})'(0)$ est presque égal à $0$. L'équation (\ref{equath2.4}) entraine également l'estimation  
$$ \left\| \nabla f \right\|_{L^\infty} \leq \sqrt{n+ \ep}.$$
 On déduit  de cette inégalité et du lemme \ref{a1l9}, l'existence d'une fonction $\tau(\ep)$
telle que pour tout $v$ dans $B_2$, 
\begin{equation}\label{a3e15}  |f(v)- \cos
d(x_0,v)| \leq \tau(\ep) 
\end{equation}
(le lemme \ref{a1l9} permet d'établir le résultat ci-dessus sur un sous-ensemble de  mesure presque égale à celle de la boule $B_2$, l'estimation sur la norme $L^{\infty}$ du gradient de la fonction $f$ permet de conclure).
En appliquant l'estimation (\ref{a3e15}) avec $v=z$, on en déduit par définition de $z$, l'existence
 d'une fonction $\eta(\ep)$ telle que
$$ d \geq \frac{\pi}{2} -\eta(\ep),$$ 
ce qui conclut.
\end{dem}

 \subsubsection{Propri\'et\'es de la vari\'et\'e ambiante} 
 Commen\c{c}ons par remarquer que sous les hypoth\`eses du th\'eor\`eme
 \ref{a3t2}, le diam\`etre de la vari\'et\'e ambiante est presque \'egal \`a
 $\pi$, pr\'ecis\'ement on a le
 \begin{lem}\label{a3l12}Sous les hypoth\`eses du th\'eor\`eme
\ref{a3t2}, il existe une fonction $\tau(\ep)$ telle que, si l'on note $d_{x_{0}}
 =\sup_{M}d(x_{0},y)$  (où $x_{0}$ est tel que $f(x_{0})=1$), alors
  $$d_{x_{0}} \geq \pi-\tau(\ep).$$ 
\end{lem} 
\begin{dem}Par le lemme \ref{a3l2} et par hypothèse sur le domaine $\Omega$, on a l'inégalité 
$$ \vol B(x_0, \frac{\pi}{2}-\eta(\ep)) \leq \frac{1}{2} \vol (M).$$
On déduit du théorème de Bishop-Gromov et par définition de $d_{x_0}$, l'inégalité 
$$ \frac{  V( \frac{\pi}{2}-\eta(\ep))}{V(d_{x_0})} \leq \frac{ \vol (B(x_0, \frac{\pi}{2}-\eta(\ep)))}{\vol B(x_0, d_{x_0})} \leq \frac{1}{2},$$
 ($V(r)$ désigne le volume d'une boule géodésique de rayon $r$ de la sphère canonique) ce qui conclut.
\end{dem}
 
 On d\'eduit en particulier de ce lemme (\`a l'aide d'une estimation de la fonction \og excess \fg \,  due \`a K. Grove et P. Petersen (\cite{grove1977}, lemme 1)  que le volume relatif de 
 toute boule de centre $x_{0}$ est presque \'egal au volume relatif d'une boule de m\^{e}me
 rayon dans la sph\`ere, par conséquent on a le
 \begin{cor} Sous les hypoth\`eses du th\'eor\`eme
\ref{a3t2}, il existe une fonction $\tau(\ep)$ telle que
 $$ \vol (\Omega \setminus B(x_{0},\frac{\pi}{2}-\eta(\ep))) \leq \tau(\ep).$$
  \end{cor}
 D'autre part, un r\'esultat de J. Cheeger et T. Colding (\cite{cheeger95}, th\'eor\`eme 5.12)
 montre que sous ces hypothèses la vari\'et\'e ambiante est proche, pour la distance de 
 Gromov-Hausdorff, d'un sinus produit tordu. 
\begin{theo}[\cite{cheeger95}]\label{a4t1}Il existe une fonction $\tau(\ep)$ telle que,
 pour tout élément $(M,g)$ de $\ens$ pour lequel il
existe $x$ dans $M$ v\'erifiant $\sup_{y \in M} d(x,y) \geq \pi - \ep$, 
l'application 
 $$\begin{array}{cccc} \phi_{x}  :&(M,g) & \longrightarrow
& ((0,\pi)\times N,d)   \\
 &z & \longmapsto & (d(x,z),p(z)) \end{array}$$
est une $\tau(\ep)$-approximation de Hausdorff, avec $N$ une sph\`ere g\'eod\'esique de $M$ 
de centre $x$ munie de sa distance intrins\`eque, $d$ la distance définie par (\ref{a4e10}) et
  $p(z)$ tel que  $d(z,p(z))=d(z,N)$.
\end{theo}

Par cons\'equent, sous les hypoth\`eses du th\'eor\`eme \ref{a3t2}, l'application 
$\phi_{x_{0}}$ est une $\tau(\ep)$-approximation de Hausdorff (lorsque $M$ est muni de la distance induite par la métrique riemmannienne) et $\phi_{x_{0}}(
B(x_{0},\frac{\pi}{2}-\eta(\ep)))= (0,\frac{\pi}{2}-\eta(\ep))\times N$. Pour 
terminer la preuve du th\'eor\`eme \ref{a3t2}, il reste \`a d\'emontrer que la boule
$B(x_{0},\frac{\pi}{2}-\eta(\ep))$ est presque convexe, c'est \`a dire que
\begin{equation}\label{a4e5}\phi_{x_{0}} :  (B(x_{0},\frac{\pi}{2}-\eta(\ep))
,d^{x_{0},(\pi-\eta(\ep))/2})
 \longrightarrow ((0,\frac{\pi}{2})\times N,d)
 \end{equation}
 est une $\tau(\ep)$-approximation.
 
 \subsubsection{Propri\'et\'e de presque convexit\'e}

 Pour d\'emontrer la propriété (\ref{a4e5}), il suffit d'\'etablir la
  \begin{prop}\label{a9p1}Sous les hypoth\`eses du th\'eor\`eme
\ref{a3t2}, il existe des fonctions $\alpha (\ep)$ et $\tau(\ep)$ telles que pour tout $u,v$ dans $B(x,\frac{\pi}{2}-\alpha (\ep)),$
  $$ d^{x_{0},(\pi-\alpha(\ep))/2}(u,v) \leq d_g(u,v) + \tau(\ep),$$
   o\`u $d_g$ est la distance induite par la m\'etrique 
riemannienne $g$ sur $M$. 
\end{prop} 

\begin{dem}
 D'apr\`es le lemme \ref{a3l2}, $\Omega$ contient la boule $B(x_0, \frac{\pi}{2}-\eta (\ep))$. 
 Nous allons montrer  que pour $\alpha (\ep)$ $(\geq \eta (\ep)) $ bien choisi, il existe une
 fonction $r (\ep)$  telle que pour tout $u,v$ dans $B(x_{0},\frac{\pi}{2}-\alpha (\ep))$,
 il existe $\tilde{u}, \tilde{v} $ v\'erifiant $d(\tilde{u},u) \leq r(\ep)$ et 
 $d(\tilde{v},v) \leq r(\ep)$ tels que
 \begin{equation}\label{a4e6}
  d^{x_{0},\frac{\pi-\alpha (\ep)}{2}}(\tilde{u},\tilde{v})= d_{g}(\tilde{u},
 \tilde{v}),
 \end{equation}
 ce qui d\'emontrera la proposition.

 Pour démontrer  l'égalité (\ref{a4e6}), nous aurons besoin du lemme suivant  (nous renvoyons à \cite{bertrand} pour une démonstration).
 \begin{lem}\label{a4l2}Il existe 
 une fonction $\tau (\ep)$ telle que
 pour tout élément $(M,g)$ de $\ens$ pour lequel il existe $ x$ appartenant à $M$
tel que $\sup_{y \in M }d(x,y) \geq  \pi -\ep$,
 il existe un vecteur $(a_{i}(x))_{1 \leq i \leq k}$ appartenant à $\mathbb{S}^{k-1}$ tel que
\begin{equation}\label{a3e17} || \cos d_x - \sum_{i=1}^{k} a_{i}(x)f_{i}||
_{L^{\infty}} \leq \tau (\ep),
\end{equation}
o\` u $d_x$ est la fonction distance à $x$, $(f_{i})$ est une famille orthogonale de fonctions propres associ\'ees aux valeurs
propres $(\lambda_{i}(M))_{1 \leq i \leq k}$, normalis\'ees par $\im f_{i}^{2}=
\frac{1}{n+1}$ et $k= \max \{i \in \mathbb{N};
\lambda_i(M) \leq n + \sqrt{\ep}\}$.
\end{lem}
\begin{rem}L'entier $k$ est supérieur ou égal à 1 d'après \cite{cheng1975}. S. Gallot a montr\'e dans \cite{gallot83} (proposition 2.4), qu'il existe
une constante positive $C(n)$ telle que pour tout élément $ (M,g)$ de $\ens $, $\lambda_{n+2}(M) \geq n+ C(n)$. Par cons\'equent, on peut supposer $ k \leq 
n+1$ dans le lemme pr\'ec\'edent.
\end{rem}
Notons $\ov[f]=\sum_{i=1}^k a_i(x)f_i$. Le lemme \ref{a4l1} implique l'existence d'une fonction $\tau (\ep)$ telle que
 \begin{equation}\label{a9e10}
  || \hes \, (\sum_{i=1}^k a_i(x)f_i) +(\sum_{i=1}^k a_i(x)f_i)g ||_{L^2} 
  \leq \tau(\ep). 
\end{equation}
Par cons\'equent, en appliquant les lemmes \ref{a1l1} et \ref{a1l9.5} \`a la fonction $ \ov[f]$, on en d\'eduit l'existence de fonctions $r(\ep)$, $R(\ep)$ et $\tau_1(\ep)$ pour lesquelles on a la propriété suivante :\\ 
Pour tout $u,v$ dans $M$ tel que $R(\ep)\leq d(u,v)\leq \pi-R(\ep),$ il existe 
 $\tilde{u},\tilde{v}$ admettant une unique g\'eod\'esique minimisante $\gamma$
les reliant et v\'erifiant $d(u,\tilde{u})<r(\ep)$,  $
d(v,\tilde{v})<r(\ep)$, tels que
 $$\left|\ov[f](\gamma(t))-\left(\cos d(x,\tilde{u})\frac{\sin (l-t)}{\sin
(l)}+ \cos d(x,\tilde{v})\frac{\sin (t)}{\sin (l)}\right)\right| \leq
\tau_1(\ep), $$ pour tout $t$ dans $[0,l]$, avec $l=d(\tilde{u},\tilde{v})$. On déduit alors de l'in\'egalit\'e (\ref{a3e17}) l'existence d'une fonction $\tau_2 (\ep)$ telle que
$$ \cos d(x,\gamma(t)) \geq \left(\cos
d(x,\tilde{u})\frac{\sin (l-t)}{\sin (l)}+ \cos
d(x,\tilde{v})\frac{\sin (t)}{\sin (l)}\right) -\tau_2 (\ep). $$ 
Puis, en remarquant que  $ \frac{\sin (l-t)}{\sin (l)}+\frac{\sin (t)}{\sin
(l)}\geq 1$, on en d\'eduit que pour tout $t$ dans $[0,l]$,
 $$ \cos d(x,\gamma(t)) \geq \mbox { min }\left(\cos d(x,\tilde{u}),\cos
d(x,\tilde{v})\right)-\tau_2(\ep),$$
ce qui permet d'\'etablir (\ref{a4e6}) et termine la  d\'emonstration.
\end{dem}
\section{Domaines g\'eod\'esiquement convexes}

Dans cette partie, nous démontrons le théorème \ref{a8t1}. Rappelons que la seconde forme fondamentale (et par
cons\'equent la courbure moyenne) du bord d'un domaine régulier géodésiquement convexe est
positive ou nulle.
Par cons\'equent, la d\'emonstration du th\'eor\`eme \ref{a8t1} se d\'eduit
du th\'eor\`eme \ref{a3t2} et de la proposition qui suit. 
\begin{prop}Il existe une fonction $R(\ep)$ telle que, pour tout
domaine r\'egulier géodésiquement convexe $\Omega$ d'un élément $(M,g)$ de $\ens$, v\'erifiant 
 $$ \frac{\vol (\Omega)}{\vol (M)} \leq \frac{1}{2} \mbox{ et } \lambda_1^D(\Omega) \leq n+ \ep,$$
 on a, en conservant les notations du th\'eor\`eme \ref{a3t2}, l'inclusion
 $$ \Omega \subset B(x_0,\frac{\pi}{2}+R(\ep)).$$
\end{prop}

\begin{dem}
On conserve les notations de la deuxième partie, $x_0$ appartient à $\Omega$ et vérifie $ f(x_0)=1$,  
$ y$ appartient à $\partial\Omega $ et vérifie $d(x_0,y)=d(x_0, \partial \Omega)$ et
$$ \frac{\pi}{2}- \eta (\ep) \leq d(x_0,y) \leq 
\frac{\pi}{2}.$$

Traitons d'abord le cas des points de $\partial \Omega$ qui sont loins de $y$. Soit
$m$ appartenant à $\partial \Omega$ et supposons que $m $ v\'erifie 
$$ \pi \geq d(m,y) \geq \pi - \mu,$$
pour $\mu$ positif petit. Dans ce cas, il y'a presque \'egalit\'e dans l'in\'egalit\'e triangulaire ci-dessous  (\cite{grove90}, lemme 1),
pr\'ecis\'ement il existe une fonction $\tau (\mu ) $ (indépendante de $x_0$) telle que
$$ d(m,x_0) + d(x_0,y) \leq d(m,y) + \tau (\mu),$$
d'où
$$ d(m,x_0) \leq \frac{\pi}{2} + \tau (\mu) + \eta(\ep).$$
Par cons\'equent, on peut dans la suite de la d\'emonstration, ne considérer que les points $m$ de $\partial \Omega$ qui vérifient
$$ d(m,y) \leq \pi -\mu$$
avec $\mu$ petit qui est défini ci-dessous.

Supposons maintenant qu'il existe un point $z$ de $\partial \Omega$ tel que
\begin{equation}\label{a4e7}
 d(x_0,z) \geq \frac{\pi}{2} + \theta (\ep)
 \end{equation}
et
$$ d(y,z) \leq \pi - \mu$$
avec $\mu =\theta (\ep) + 6 r(\ep)$ 
 (où $r(\ep)$ est tel que $
\frac{\ep}{V(r(\ep))} \leq \ep^{1/2}$). \\

Nous allons montrer que pour $\theta (\ep)$ convenable, on obtient une contradiction. Pour cela, on montre qu'un point $P$ situé à mi-distance entre $z$ et $y$ contient un voisinage ouvert dont on peut minorer le volume (\ref{bidule}). On peut ensuite grâce aux résultats de la première partie, estimer la valeur de $f$ en $P$, en fonction de la distance $d(x_0,P)$. On obtient la contradiction souhaitée en estimant $d(x_0,P)$ à l'aide de la $\tau (\ep)$-approximation de Hausdorff $\phi_{x_0}$. 

Notons $z'$ un point appartenant à une géodésique minimisante reliant $x_0$ à $z$ et tel que $d(z,z')=r(\ep)$,
définissons de manière analogue $y'$ vérifiant $ d(y,y')=4r(\ep)$. Notons $d_1= d(y',z')$. On déduit de l'inégalité triangulaire
$$d_1 \geq d(x_0,z) - d(x_0,y) + 3r(\ep),$$
 ce qui entraine par (\ref{a4e7}) et comme $d(x_0,y) \leq \pi /2$,
 $$ d_1 \geq \theta (\ep)+ 3r(\ep).$$
 
 Soit deux réels $s,S$ positifs, $p$ appartenant à $M$ et $\Gamma$ une partie mesurable de $S_p(M)$ la
sphère unitaire tangente en $p$. On note
$$A_{s,S}^{\Gamma}(p)=\{\gamma_{\vec{u}}(t); \vec{u} \in \Gamma \mbox{ et } t \in [s,S]\}$$ 
et ${\ov[A]}_{s,S}^{\Gamma}$ l'ensemble correspondant sur la sph\`ere. 
 
 A l'aide de l'inégalité triangulaire, on montre 
 \begin{equation}\label{a9e1} B(y',r(\ep)) \subset
A_{d_1-r(\ep),d_1+r(\ep)}^{\Gamma}(z') \subset B(y',3r(\ep)),
\end{equation} avec $\Gamma$ la trace sur $S_{z'}(\Omega)$ de
$B(y',r(\ep))$. Par construction, la boule $B(y',3r(\ep))$ est contenue dans $\Omega$. Par conséquent, en utilisant la  convexité de $\Omega$, on en déduit 
$$A_{\frac{d_1}{2}-r(\ep),\frac{d_1}{2}+r(\ep)}^{\Gamma}(z') \subset
\Omega.$$

Estimons maintenant le volume de $A_{\frac{d_1}{2}-r(\ep),\frac{d_1}{2}+r(\ep)}^{\Gamma}(z')$. Soit $s,S,r,R$ des r\'eels v\'erifiant $0<s<S$, $0<r<R$, $s<r$, $S<R$ et $p$ dans $M$. Le théorème de Bishop-Gromov implique
 $$ \frac{\vol(A_{s,S}^{\Gamma}(p))}{\vol (A_{r,R}^{\Gamma}(p))} \geq \frac{\vol
({\ov[A]}_{s,S}^{\Gamma})}{\vol ({\ov[A]}_{r,R}^{\Gamma})}.$$
Donc, on a en particulier
\begin{equation}\label{a9e3}
\vol (A_{\frac{d_1}{2}-r(\ep),\frac{d_1}{2}+r(\ep)}^{\Gamma}(z')) \geq
\vol (A_{d_1-r(\ep),d_1+r(\ep)}^{\Gamma}(z'))\times
\frac{\vol({\ov[A]}_{\frac{d_1}{2}-r(\ep),\frac{d_1}{2}+r(\ep)}^{\Gamma
})}{\vol({\ov[A]}_{d_1-r(\ep),d_1+r(\ep)}^{\Gamma})}. 
\end{equation}
Sur la sphère canonique, on a l'égalité
 $$ 
\frac{\vol ({\ov[A]}_{\frac{d_1}{2}-r(\ep),\frac{d_1}{2}+r(\ep)}^{\Gamma})}{{
\vol (\ov[A]}_{d_1-r(\ep),d_1+r(\ep)}^{\Gamma})}= \frac{
\int_{\frac{d_1}{2}-r(\ep)}^{\frac{d_1}{2}+r(\ep)}\sin^{n-1} (t)dt}{
\int_{d_1-r(\ep)}^{d_1+r(\ep)}\sin^{n-1} (t)dt}.$$ 
Par cons\'equent, en remarquant que pour tout $t$ dans $[0,\pi]$, 
\begin{equation}\label{a9e6} \frac{\sin^{n-1} (\frac{t}{2})}{\sin^{n-1}
(t)} \geq \frac{1}{2^{n-1}}, 
\end{equation}
 on en d\'eduit en intégrant
 $$ 
\frac{\vol({\ov[A]}_{\frac{d_1}{2}-r(\ep),\frac{d_1}{2}+r(\ep)}^{\Gamma})}{\vol ({
\ov[A]}_{d_1-r(\ep),d_1+r(\ep)}^{\Gamma})} \geq \frac{1}{2^n}.$$
D'où, par (\ref{a9e1}) et (\ref{a9e3})
\begin{equation}\label{bidule}
  \vol (A_{\frac{d_1}{2}-r(\ep),\frac{d_1}{2}+r(\ep)}^{\Gamma}(z')) \geq 
 \frac{1}{2^n} \vol (B(y',r(\ep))).
 \end{equation}
Par choix de $r(\ep)$, on déduit des lemmes \ref{a4l1} et \ref{a1l1}
 (en notant $A=\\
 A_{\frac{d_1}{2}-r(\ep),\frac{d_1}{2}+r(\ep)}^{\Gamma}(z')$)
 
\begin{multline*}
 \frac{1}{\vol(A \times B(x_0,r(\ep)))} \int_{ A \times B(x_0,r(\ep)}
 \left(\int_{0}^{1} |(f\circ \gamma_{uv})''(t) +f( \gamma_{uv}(t))|^2dt
 \right)dudv \\
 \leq \tilde{C}(n)\ep^{\frac{1}{2}}. 
 \end{multline*}
Par cons\'equent, comme $x_0$ est un point où $f$ réalise son maximum, on en 
déduit, comme dans la preuve du lemme \ref{a3l2}, qu'il existe une fonction 
$\tau'(\ep)$ telle que, pour tout $v$ appartenant à $A$,
\begin{equation}\label{a4e8}
\left|f(v)- \cos d(x_0,v)\right| \leq \tau' (\ep).
\end{equation}

 Nous allons maintenant estimer $\cos d(x_0,v)$ pour $v$ appartenant à
  $A$ fixé, en utilisant la $\tau (\ep)$-approximation de Hausdorff $\phi_{x_0}$ du théorème  \ref{a4t1}. Notons $\ov[x],\ov[y],\ov[z]$ et $\ov[v]$ l'image des
points $x_0,y',z'$ et $v$ par $\phi_{x_0}$, ces points v\'erifient
\`a $\tau(\ep)$ pr\`es les m\^emes relations métriques que les points $x,y',z'$ et $v$. Par définition de l'application
$\phi_{x_0}$ et de la distance sur le produit tordu $ (0,\pi)\times N$,
on peut, pour  estimer la distance de $\ov[x]$ \`a $\ov[v]$, supposer que l'on est sur la sph\`ere canonique
 (voir la remarque 1.47, page 196 de \cite{cheeger95}).

 Afin d'estimer la distance $d(\ov[x],\ov[v])$ sur la sphère canonique, on
suppose dor\'enavant que les quantités $\frac{\tau(\ep)}{\theta(\ep)},\frac{\tau'(\ep)}{\theta(\ep)},\frac{\eta(\ep)}{\theta(\ep)},\frac{r(\ep)}{\theta(\ep)}$ tendent vers $0$ quand $\ep$ tend vers $0$, de sorte que les termes $\tau (\ep),\tau'(\ep), \eta(\ep)$ et $r(\ep)$ sont n\'egligeables devant $\theta(\ep)$, cela assure également que $\cos (\frac{d(\ov[z],\ov[y])}{2})$ n'est pas trop petit, puisque, à des termes négligeables près, on a $  d(\ov[z],\ov[y]) \leq \pi - \theta(\ep)$. Sur
la sph\`ere canonique et en utilisant l'hypothèse (\ref{a4e7}), un calcul montre qu'il existe un r\'eel $C >0$ tel que 
$$ d(\ov[x],\ov[v]) \geq \frac{\pi}{2} + C \theta(\ep).$$
Par cons\'equent
$$ d(x,v) \geq  \frac{\pi}{2} + C \theta(\ep)-\tau(\ep).$$
\`A l'aide de (\ref{a4e8}), on en d\'eduit 
 $$ f(v) \leq \cos \, (\frac{\pi}{2} + C
\theta(\ep)-\tau(\ep)) +\tau'(\ep)$$ 
et donc quitte \`a supposer $\theta(\ep)$ assez grand  
 $$ f(v) <0,$$
 ce qui est absurde.
\end{dem}

\bibliographystyle{plain}

\end{document}